\documentclass{birkjour}
\usepackage [latin1]{inputenc}

\newtheorem{thm}{Theorem}[section]

\newtheorem{lem}[thm]{Lemma}

\theoremstyle{definition}

\theoremstyle{remark}

\numberwithin{equation}{section}

\newcommand{\BibTeX}{B\kern-0.1emi\kern-0.017emb\kern-0.15em\TeX}
\newcommand{\XYpic}{$\mathrm{X\kern-0.3em\raisebox{-0.18em}{Y}}$-$\mathrm{pic}\,$}

\newcommand{\cl}{C \kern -0.1em \ell}  

\newcommand{\ed}{\end{document}}
\begin{document}
\title[A note on Ricci soliton in almost Kenmotsu manifold]
{A note on Ricci soliton in almost Kenmotsu manifold}

\author[Fatemah Mofarreh and U. C. De]
{Fatemah Mofarreh and Uday Chand De}

\address
{ Mathematical Science Department, Faculty of Science, Princess Nourah bint Abdulrahman University, Riyadh, 11546, Kingdom of  Saudi Arabia .}

\address
{Department of Pure Mathematics, University of Calcutta
35, Ballygaunge Circular Road
Kolkata 700019, West Bengal, India.}
\email {uc$_{-}$de@yahoo.com}

\footnotetext {$\bf{2010\ Mathematics\ Subject\ Classification\:}.$ 53C25, 53D15.
\\ {Key words and phrases: Almost Kenmotsu manifold, Ricci soliton.\\
}}

\begin{abstract}
In this short note it is established that there does not exist Ricci soliton with the Reeb potential vector field in an almost Kenmotsu manifold (briefly, $\mathcal{AKM}$).
\end{abstract}
\maketitle
\numberwithin{equation}{section}
\newtheorem{theorem}{Theorem}[section]
\newtheorem{lemma}[theorem]{Lemma}
\newtheorem{proposition}[theorem]{Proposition}
\newtheorem{corollary}[theorem]{Corollary}
\newtheorem*{remark}{Remark}
\newtheorem{Agreement}[theorem]{Agreement}
\newtheorem{definition}[theorem]{Definition}
\newtheorem{example}[theorem]{Example}

\section{\textbf{Introduction}}

A Ricci soliton, a natural generalization of an Einstein metrics in a Riemannian manifold $(\mathcal{M},g)$ is a triplet $(g, \mathcal{V},\lambda)$, with $g$, a Riemannian metric, $\mathcal{V}$ a smooth vector field (called the potential vector field) and $\lambda$ a constant such that
\begin{equation}\label{100}
\pounds_{\mathcal{V}}g+2S-2\lambda g=0,
\end{equation}
where $\pounds_{\mathcal{V}}g$ and $S$ denotes the Lie derivative of $g$ along a vector field $\mathcal{V}$ and the Ricci tensor of type $(0,2)$. The above-stated soliton is an Einstein metric if $\mathcal{V}$ is zero or Killing. The Ricci soliton is termed shrinking, steady, or expanding according as $\lambda$ is positive, zero, or negative, respectively. Compact Ricci solitons are nothing but the fixed points of the Ricci flow $\frac{\partial}{\partial t}g=-2S$ projected from the space of metrics onto its quotient modulo diffeomorphisms and scalings. On compact manifolds, they sometimes arise as blow-up limits for the Ricci flow. Metrics fulfilling (\ref{100}) are exceptionally helpful in modern physics. Theoretical physicists have been investigating the condition of Ricci soliton in connection with string theory. On a compact manifold this soliton has constant curvature in 2-dimension(Hamilton, \cite{hn}) and also in 3-dimension(Ivey, \cite{iy}).\par

Recently, Cho \cite{cj1} obtained some interesting results about Ricci solitons in almost contact and contact geometry.  Also, Ricci solitons have been studied by several authors such as  Chow and Knopf \cite{ck}, De et al. \cite{dsm17}, De and Mandal \cite{de4}, Deshmukh et al. \cite{desh2},  Wang and Liu \cite{wlt} and many others.\par\par

In \cite{cho}, Cho proved that there is no Ricci soliton with the Reeb potential vector field in a Kenmotsu manifold.\par
The purpose of this short note is to generalize the above result of Cho in $\mathcal{AKM}$s.

\section{\textbf{Almost Kenmotsu manifolds}}

A differentiable manifold $\mathcal{M}^{(2n+1)}$ is said to have a $(\phi, \xi, \eta)$-structure or an almost contact structure, if it admits a $(1,1)$ tensor field $\phi$, a characteristic vector field $\xi$ and a 1-form $\eta$ fulfilling \cite{bl}
\begin{equation}\label{g1}
\phi^{2}=-I+\eta\otimes\xi,\; \eta(\xi)=1, \end{equation}
 where $I$ indicates the identity endomorphism. Here, also $\phi\xi=0$ and $\eta\circ\phi=0$; both can be derived from the equation $(\ref{g1})$ easily. \par
If a manifold $\mathcal{M}$ with a $(\phi, \xi, \eta)$-structure admits a Riemannian metric $g$ such that
\begin{equation}\nonumber g(\phi E,\phi F)=g(E,F)-\eta(E)\eta(F), \end{equation}
for any vector fields $E$, $F$ of $T_{p}\mathcal{M}^{2n+1}$, then $\mathcal{M}$ is said to be an almost contact metric manifold. The fundamental 2-form $\Phi$ on an almost contact metric manifold is
defined by $\Phi(E,F)=g(E,\Phi F)$ for any $E$, $F$ of $T_{p}\mathcal{M}^{2n+1}$. The condition for an almost contact metric manifold being normal is equivalent to the vanishing of the $(1,2)$-type torsion tensor $\mathcal{N}_{\phi}$, defined by $\mathcal{N}_{\phi}=[\phi,\phi]+2d\eta\otimes\xi$, where $[\phi,\phi]$ is the Nijenhuis torsion of $\phi$. \cite{bl}.
Recently, in (\cite{dp},\cite{dp2}), almost contact metric manifolds such that $\eta$ is closed and $d\Phi=2\eta\wedge\Phi$ are investigated and they are named $\mathcal{AKM}$s. Obviously, a normal $\mathcal{AKM}$ is a Kenmotsu manifold. Also Kenmotsu manifolds can be characterized by
\begin{equation}\nonumber (\nabla_{E}\phi)F=g(\phi E,F)\xi-\eta(F)\phi E, \end{equation}
for any vector fields $E$, $F$. It is well known \cite{kk} that a Kenmotsu manifold $\mathcal{M}^{2n+1}$ is locally a warped product $I\times_{f}\mathcal{N}^{2n}$ where $\mathcal{N}^{2n}$ is a K\"{a}hler manifold, $I$ is an open interval with coordinate $t$ and the warping function $f$, defined by $f=ce^{t}$ for some positive constant c. Let us denote the distribution orthogonal to $\xi$ by $\mathcal{D}$ and defined by $\mathcal{D}=Ker(\eta)=Im(\phi)$. In an $\mathcal{AKM}$, since $\eta$ is closed, $\mathcal{D}$ is an integrable distribution. Let $\mathcal{M}^{2n+1}$ be an $\mathcal{AKM}$. We denote by $h=\frac{1}{2}\pounds_{\xi}\phi$ and $l=R(\cdot, \xi)\xi$ on $\mathcal{M}^{2n+1}$. The tensor fields $l$ and $h$ are symmetric operators and obey the subsequent relations :
\begin{equation}\label{b1}
h\xi=0,\;l\xi=0,\;tr(h)=0,\;tr(h\phi)=0,\;h\phi+\phi h=0,
\end{equation}
\begin{equation}\label{b2}
 \nabla_{E}\xi=-\phi^2E-\phi hE(\Rightarrow \nabla_{\xi}\xi=0),
\end{equation}
\begin{equation}\label{b3} \phi l\phi-l=2(h^{2}-\phi^{2}),
\end{equation}
\begin{equation}\label{b4}
R(E,F)\xi=\eta(E)(F-\phi hF)-\eta(F)(E-\phi hE)+(\nabla_{F}\phi h)E-(\nabla_{E}\phi h)F,
\end{equation}
for any vector fields $E$, $F$. The $(1,1)$-type symmetric tensor field $h'=h\circ\phi$ is anticommuting with $\phi$ and $h'\xi=0$. Also it is clear that \cite{dp}
\begin{equation}\label{b5}
 h=0\Leftrightarrow h'=0,\;\;h'^{2}=(k+1)\phi^2(\Leftrightarrow h^{2}=(k+1)\phi^2).
 \end{equation}

$\mathcal{AKM}$s have been studied by  several authors such as  Dileo and Pastore (\cite{dp}, \cite{dp2}), Wang and Liu (\cite{wa}, \cite{wlt}), De and Mandal \cite{de4} and many others.\par

Here we state a lemma which will be used later.

\begin{lem}\cite{cj1}
If $(g,\mathcal{V})$ is a Ricci soliton of a Riemannian manifold, then we acquire
\begin{eqnarray}\label{l1}
\frac{1}{2}||\pounds_{\mathcal{V}}g||^{2}=dr(\mathcal{V})+2div(\lambda \mathcal{V}-Q\mathcal{V}),
\end{eqnarray} where $r$ denotes the scalar curvature and $Q$ denotes the Ricci operator defined by $S(E,F)=g(QE,F)$.
\end{lem}

\section{\textbf{Ricci solitons}}

In this section we characterize $\mathcal{AKM}$ admitting a Ricci soliton. Now suppose that an $\mathcal{AKM}$ admits a Ricci soliton. Then from $(\ref{100})$, we obtain
\begin{eqnarray}\label{1}
\frac{1}{2}[g(\nabla_{E}\xi, F)+g(\nabla_{F}\xi, E)]+S(E,F)-\lambda g(E,F)=0.\end{eqnarray}
Utilizing $(\ref{b2})$  and $(\ref{g1})$ in $(\ref{1})$ yields
\begin{eqnarray}\label{2}
&&
\frac{1}{2}[g(E-\eta(E)\xi-\phi hE, F)+g(F-\eta(F)\xi-\phi hF,E)]\nonumber\\&&+S(E, F)-\lambda g(E,F)=0.
\end{eqnarray}
Using $(\ref{b1})$ in $(\ref{2})$ implies
\begin{eqnarray}\nonumber
\frac{1}{2}[2g(E,F)-2\eta(E)\eta(F)-2g(\phi hE, F)]+S(E,F)-\lambda g(E,F)=0,
\end{eqnarray}
that is,
\begin{eqnarray}\label{3}
g(E,F)-\eta(E)\eta(F)-g(\phi hE, F)+S(E,F)-\lambda g(E,F)=0.
\end{eqnarray}
From the preceding equation it follows that
\begin{eqnarray}\label{4}
QE=(\lambda-1)E+\eta(E)\xi+\phi hE.
\end{eqnarray}
Replacing $E$ by $\xi$ in $(\ref{4})$ we infer
\begin{eqnarray}\label{5}
Q \xi=\lambda \xi.
\end{eqnarray}
Let $\{e_{1},e_{2},e_{3},...,e_{2n+1}\}$ be a local orthonormal basis of the tangent space at a point of the manifold $\mathcal{M}$. Then by putting $E=F=e_{i}$ in $(\ref{3})$ and taking summation over $i$, 
$1\leq i  \leq (2n+1)$, we have
\begin{equation}\label{200}
r=\lambda(2n+1)-2n.
\end{equation}
Therefore the scalar curvature is constant, as $\lambda$ is constant.\\
Since
\begin{eqnarray}
\label{6}
\frac{1}{2}||\pounds_{\mathcal{V}}g||^{2}=dr(\mathcal{V})+2 div (\lambda \mathcal{V}-Q\mathcal{V}),
\end{eqnarray}
so replacing $\mathcal{V}$ by $\xi$ and in view of the fact the scalar curvature $r$ is constant, we get
\begin{eqnarray}
\label{7}
\frac{1}{2}||\pounds_{\xi}g||^{2}=2 div (\lambda \xi-Q\xi).
\end{eqnarray}
Making use of $(\ref{5})$ in $(\ref{7})$ we obtain
\begin{equation}\nonumber
||\pounds_{\xi}g||^{2}=0,
\end{equation}
which implies that $\xi$ is a Killing vector field. But in an $\mathcal{AKM}$ $\xi$ can never be a Killing vector field.

This leads to the following:
\begin{theorem}
There does not exist Ricci soliton in an almost Kenmotsu manifold.
\end{theorem}

\end{document}